\documentclass[twoside,reqno,10pt]{article}
\usepackage{tikz}

\usepackage{amsopn,amsfonts,amscd,amssymb,latexsym,amsmath,amsthm}
\usepackage{hyperref}

\newtheorem{thm}{Theorem}[section]

\newtheorem{lem}[thm]{Lemma}

\newtheorem{rem}[thm]{Remark}
\newtheorem{exmp}[thm]{Example}

\theoremstyle{definition}

\numberwithin{equation}{section}

\begin{document}
\begin{center}
{\Large On a Carleman formula for lunes}
\end{center}
\begin{center}
{\large D. Fedchenko\footnote{Institute of Mathematics, Siberian Federal University, Svobodny Prospect 79, Krasnoyarsk, 660041, Russia. E-mail: dfedchenk@gmail.com.\\
The research of the author was done in the framework of the Mikhail Lomonosov Fellowship which is supported by the Russian Ministry of Education and the Deutsche Forschungsgemeinschaft.
}}
\end{center}

\begin{abstract}
In this paper we consider a simple formula for analytic continuation in a domain $D
\subset \mathbb{C}$ of special form.
\end{abstract}

\section*{Introduction}
Integral representations of holomorphic functions solve the
classical problem on restoring a holomorphic function in a
domain $D$ by its values on $\partial D$. Connected with this
problem, there is another one: restore holomorphic function in $D$
by its values on a set $\Gamma \subset \partial D$.

The first result of such type was obtained by Carleman \cite{C} for a plain domain $D$ of special form. His idea of using {\it extinguishing} functions
was developed in the article of Goluzin and Krylov \cite{GK} and by Fok and Kuni \cite{FK} for simply connected plain domains of special form. This
method provides an {\it extinguishing} function for any subset $\Gamma$ of $\partial D$ of positive measure. Another method was proposed by Lavrent'ev in 1956 \cite{L}.

All these formulas and some their applications can be found in the book of Aizenberg \cite{A}.

In this notice we show a simple trick to construct
Carleman formula for lunes.

\section{Analytic continuation}

Let $\mathbb{C}$ be the complex plane of $z = x + \imath y$, where $x,y \in \mathbb{R}$. For an open set $D$, denote by $\mathcal{O}(D)$ the space of functions holomorphic in $D$. Write $\overline \partial$ for the Cauchy-Riemann operator
$$
\overline \partial = \frac{1}{2} \left( \frac{\partial}{\partial x} + \imath \frac{\partial}{\partial y} \right)
$$
 in $\mathbb{C}$.

Let $\Gamma$ be a subset of $\partial D$. The problem of analytic continuation from $\Gamma$ into $D$ consists in the following. Given a function $u_0 \in C(\Gamma)$,  find a function $u \in \mathcal{O}(D) \cap C(D \cup \Gamma)$ such that

$$
\left\{ \begin{array}{l}
\overline \partial u = 0_{\phantom{0}} \mbox{ \it in } D, \\
\phantom{\partial} u = u_0 \mbox{ \it on }
\Gamma.
\end{array} \right.
$$
One easily specifies this problem within Cauchy problems for solutions of elliptic equations.

\section{Carleman formula}

Let $\Gamma$ be a smooth curve in the unit disk $B(0,1)$ dividing $B(0,1)$ into two domains. Denote by $D$ those of these domains which does not contain the origin. Such domains $D$ are referred to as lunes. And let
\begin{align*}
\gamma(t) &= (x(t), y(t)), \, t \geq 0, \\
\gamma(0) &= 0
\end{align*}
be a smooth curve with
end point the origin, which lies in $B(0,1) \setminus D$ (see
Fig.~1). Let $\Gamma$ and $\gamma$ intersect transversally at the
origin. Fix now any curve $\gamma$ with these properties.

\begin{center}

\begin{tikzpicture} \draw[step=.5cm,gray,very thin]
(-1.5,-1.5) grid (1.5,1.5);

\draw[line width=1pt] (0cm,-1.5cm) arc (-90:90:1.5cm) node[below=5pt,fill=white] {$D$} (0,1.5cm); 

\draw[dashed] (0cm,1.5cm) arc (90:270:1.5cm);%
\draw (0cm,0cm) arc (90:270:0.25cm);%
\draw (0cm,-0.5cm) arc (90:-90:0.25cm);%
\draw (0,-0.5cm) circle (1cm);%

\draw[color=red, line width=1pt] (0cm,1.5cm) arc (90:270:0.75cm)
-- node[right=10pt] {$\Gamma$} (0,0pt);%
\draw[color=red, line width=1pt] (0cm,-1.5cm) arc (-90:90:0.75cm);%
\draw[help lines] (0,-0.5cm) circle (1cm);%
\filldraw [gray] (0,0) circle (1pt) node[above=1pt, fill=white] {$0$}(0pt,0pt);%
\filldraw [gray] (0,-0.5cm) circle (1pt); \draw (0cm,-1cm) arc
(90:130:0.7cm) node[left=0pt, fill=white]
{$\gamma(t)$}(0pt,1.5cm);%


\end{tikzpicture}

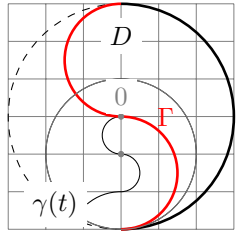
Fig. 1. Example of $D$
\end{center}

Denote by $B(\gamma(t))$ the circle in $\mathbb{C}$ with
center at the point $\gamma(t)$ and with radius $\textrm{dist}
(\gamma(t), \partial B(0,1))$.

\begin{lem}\label{lemma_a_exist}
For any $z\in D$, there exists point $\gamma(t_0)$ such that
$z \in B(\gamma(t_0))$.
\end{lem}
\textit{Proof.} Proof by contradiction: suppose there is a point $z \in D$ such that $z \notin B(\gamma(t))$ for all $t \geq 0$, then
$$
\left \{ \begin{array}{l} | z - \gamma(t) | \geq
\textrm{dist}(\gamma(t), \partial B(0,1)), \mbox{ for all } t \geq 0; \\
| z | < 1.
\end{array} \right.
$$
Let $t \rightarrow 0$. Then we have
$$
\left \{ \begin{array}{l}%
| z| \geq 1, \\
| z | < 1.
\end{array} \right.
$$
Obtained contradiction proofs this lemma. \hfill $\square$

Take any numerical sequence $\gamma(t_N) = (x(t_N), y(t_N))$, $N
\in \mathbb{N}$, $t_N \geq 0$ such that $\lim \limits_{N
\rightarrow \infty} \gamma(t_N) = 0$.

\begin{thm}\label{theorem.1}
If $u \in \mathcal{O}(D) \cap C(\overline D)$, then the formula
$$
u(z) = \frac{1}{2 \pi \imath} \lim_{N \rightarrow \infty}
\int_\Gamma \frac{u(\zeta)}{\zeta - z} \left( \frac{z -
\gamma(t_N)}{\zeta - \gamma(t_N)} \right)^{N+1} d\zeta,
$$
holds for any point $z \in D$, where the convergence is uniform in $z$ on compact subsets of $D$.
\end{thm}

{\it Proof.} Fix $z \in D$. Choose $N$ large enough, so that $z \in B(\gamma(t_N))$. Expand the Cauchy kernel as Laurent series in the variable $\zeta$ in the complement of the disk $B(\gamma(t_N))$ by
$$
\frac{1}{\zeta - z} = \sum_{k=0}^\infty \frac{(z - \gamma(t_N))^{k \phantom{+ 1}}}{(\zeta - \gamma(t_N))^{k+1}}.
$$

Consider the sequence of kernels
\begin{equation}\label{Carleman_Kernels}
\mathfrak{C}_N(\zeta,z) = \frac{1}{2 \pi \imath} \left(
\frac{1}{\zeta - z} - \sum_{k=0}^N \frac{(z -
\gamma(t_N))^{k \phantom{+ 1}}}{(\zeta - \gamma(t_N))^{k+1}} \right)
\end{equation}
which we call Carleman kernel.

Using the geometric sum formula, we make the following
transformation of Carleman kernels (\ref{Carleman_Kernels})
$$
\mathfrak{C}_N(\zeta,z) = \frac{1}{2 \pi \imath} \frac{1}{\zeta -
z} \left( \frac{z - \gamma(t_N)}{\zeta - \gamma(t_N)}
\right)^{N+1}.
$$

Let $N \rightarrow \infty$. Write
\begin{multline}\label{u(z).formula}
\frac{1}{2 \pi \imath} \lim_{N \rightarrow \infty} \left[
\int_{\partial D} \frac{u(\zeta)}{\zeta - z} \left( \frac{z -
\gamma(t_N)}{\zeta - \gamma(t_N)} \right)^{N+1} d\zeta \,\,\, - \right. \\
\left. \int_{\partial
D \setminus \Gamma} \frac{u(\zeta)}{\zeta - z} \left( \frac{z -
\gamma(t_N)}{\zeta - \gamma(t_N)} \right)^{N+1} d\zeta \right] = \\
\frac{1}{2 \pi \imath} \lim_{N \rightarrow \infty}
\int_{\Gamma} \frac{u(\zeta)}{\zeta - z} \left( \frac{z -
\gamma(t_N)}{\zeta - \gamma(t_N)} \right)^{N+1} d\zeta.
\end{multline}

Fix any compact $K \subset D$. We see that
the second integral from the left-hand side in (\ref{u(z).formula}) tends to zero uniformly on $K$, because $q_K =
\max \limits_{z\in K, \zeta \in \partial D \setminus \Gamma}
\frac{|z|}{|\zeta|}<1$. And from the holomorphy of the function
$u(\zeta) \left( \frac{z - \gamma(t_N)}{\zeta - \gamma(t_N)}
\right)^{N+1}$, for all $N \in \mathbb{N}$, and form the Cauchy formula we conclude that the integral over $\partial D$ equals $u(z)$ in the domain $D$. 

We thus arrive at the desired formula
$$
u(z) = \frac{1}{2 \pi \imath} \lim_{N \rightarrow \infty}
\int_{\Gamma} \frac{u(\zeta)}{\zeta - z} \left( \frac{z -
\gamma(t_N)}{\zeta - \gamma(t_N)} \right)^{N+1} d\zeta.
$$
\hfill$\square$

\begin{exmp}
Let $D=\{ |z|<1, \Re z > 0 \}$ and $\Gamma = \{ |z|<1, \Re z = 0
\}$. As $\gamma(t_N)$ we take the sequence
$(-1/N,0)$, $N\in \mathbb{N}$. If $u \in \mathcal{O}(D) \cap C(\overline D)$, then for any point $z \in D \cup \Gamma$ the formula
$$
u(z) = \frac{1}{2 \pi \imath} \lim_{N \rightarrow \infty}
\int_\Gamma \frac{u(\zeta)}{\zeta - z} \left( \frac{Nz + 1}{N\zeta + 1} \right)^{N+1} d\zeta
$$
holds, where the integral converges uniformly on compact subsets of $D$.
\end{exmp}

\begin{rem}
If $z=0$ does not belong to $\overline D$ then our Carleman formula recovers one of the formulas in \cite{A}.
\end{rem}


\begin{thebibliography}{XXXX}

\bibitem[C]{C} T. Carleman, \textit{Les fonctions quasianalytiques}, Paris: Gauthier-Villars. (1926).

\bibitem[GK]{GK} G. Goluzin and V. Krylov, \textit{Generalized Carleman formula and its application to analytic continuation of functions} Mat. Sb., {\bf 40}, 144~-~149 (1933).

\bibitem[FK]{FK} V. Fok and F. Kuni, {\it On the cutting function in dispersion relations}, Dokl. Akad. Nauk SSSR {\bf 127} (1959), 1195-1198 (Russian)

\bibitem[LRS]{L}  M. M. Lavrent'ev, V. G. Romanov and S. P. Shishatskii. \textit{Ill-Posed Problems of Mathematical Physics and Analysis} [in Russian], Nauka, Moscow (1980); English transl.: Amer. Mathem. Soc., Providence Vol. 64 (1986).

\bibitem[A]{A}  L. Aizenberg, \textit{Carleman's formulas in complex analysis}, Kluwer Acad. Publ. (1993).
\end{thebibliography}
\end{document}